# POISSON LIMITS OF SUMS OF POINT PROCESSES AND A PARTICLE-SURVIVOR MODEL

By Matthew O. Jones and Richard F. Serfozo

*Austin Peay State University and Georgia Institute of Technology*

We present sufficient conditions for sums of dependent point processes to converge in distribution to a Poisson process. This extends the classical result of Grigelionis [*Theory Probab. Appl.* **8** (1963) 172–182] for sums of uniformly null point processes that have Poisson limits. Included is an application in which a particle-survivor point process converges to a Poisson process. This result sheds light on the "surprising" Poisson limit of the species competition process of Durrett and Limic [*Stochastic Process. Appl.* **102** (2002) 301–309].

**1. Introduction.** Tilman [7] studied a model for species coexistence where species with higher death rates are superior competitors (as can be the case with parasites that kill their hosts quickly). He described the dynamics of such systems with differential equations and provided solutions describing the relative abundances of the species in equilibrium. Interestingly, he was able to show that inferior competitors with suitable traits can endure, regardless of the number of species in the system.

To further examine the number of species that can coexist, May and Nowak [5] considered a model where new species continually arrive to a system, possibly due to mutations or arrivals of new species from outside the system. They assumed sufficient time would elapse between arrivals for the system to come to equilibrium. Among other things, the authors were interested in the limiting distribution of the numbers of different species with death rates in subintervals of $(0, 1)$. In order to study this, they introduced a toy model where an arriving species with death rate $v$ immediately annihilates any species with a smaller death rate with probability $\alpha$, independently of everything.

Durrett and Limic [1] gave a formal treatment of this model. They studied a system in which particles arrive to the interval $(0,1)$ according to a Poisson process, and the location of each particle in the interval has a uniform distribution, independent of everything else. The location of a particle represents the particle's rank, which affects how long it will survive. Whenever a new particle arrives, each particle currently in the interval with a rank below that of the arrival is independently deleted with a fixed probability. The evolution over time of the point process of particles in the interval is represented by a Markov jump process on a space of counting measures. The authors proved this Markov process converges to a Poisson process on $(0,1)$; its intensity is increasing and approaches infinity at 1.

In trying to determine why Poisson processes arise as limits in such applications, we studied several types of sums of point processes with dependencies that are abstractions of those in the particle process. This led us to a general Poisson convergence theorem for sums of dependent point processes (Theorem 1). This result is similar in structure to Grigelionis' [2] theorem that gives necessary and sufficient conditions for independent uniformly null processes to converge to a Poisson process.

The rest of this study shows that Theorem 1 is a natural framework for establishing Poisson limits of particle-survivor processes. We apply it to a system in which particles enter $\Re$ at arbitrary random times, and their locations in $\Re$ are independent and identically distributed with a general, continuous distribution. Each new arrival may delete any particle located below it with a location-dependent probability $a(x)$ (see Section 3). We show that the point process of particles converges in distribution to a Poisson process. Our results indicate that, by appropriately choosing the function $a(x)$, one can allow an arbitrarily large number of particles to exist in any subinterval of $(0,1)$, on the average. This may appear contrary to the "limiting similarity" of Tilman's model, but is actually a simple consequence of the toy models of [1, 5] since each particle has a finite lifetime w.p.1.

The contents of the remaining sections are as follows. Section 2 contains our main result for sums of point processes and several corollaries for special cases. A preliminary lemma on the convergence in distribution of products of random variables is of interest by itself. Section 3 describes the particle-survivor process described above and its Poisson limit. Sojourn times of particles in this system and properties of departures are the topics of Section 4. We end in Section 5 by giving necessary and sufficient conditions for the stationary distribution of a Markovian particle-survivor process to be that of a Poisson process.

**2. Sums of dependent point processes.** Before stating our main result, we review some terminology on point process as in [3, 4]. Let $(\Omega, \mathcal{F}, P)$ be a probability space and $\mathbb{E}$ be a Polish space endowed with its Borel $\sigma$-field $\mathcal{E}$. A *point process* $\xi$ on $\mathbb{E}$ is a mapping from $\Omega$ to the set of integer-valued measures of the form $\sum_k \delta_{x_k}(\cdot)$ that are finite on compact sets, where $\delta_x(\cdot)$ is the Dirac measure with a unit atom at $x$. This space of measures is endowed with the Borel $\sigma$-field generated by the topology of vague convergence. Then $\xi(B) = \sum_k \delta_{X_k}(B)$ represents the number of points in $B \in \mathcal{E}$, where $X_k$ are the point locations. The mean measure of $\xi$ is $\mu(B) \equiv E[\xi(B)]$, $B \in \mathcal{E}$.

The point process $\xi$ is a *Poisson process* with mean measure $\mu$ if $\xi(B_1), \ldots, \xi(B_m)$ are independent for disjoint $B_1, \ldots, B_m \in \hat{\mathcal{E}}$, and $\xi(B)$ is a Poisson random variable with mean $\mu(B)$, where $\mu$ is finite on compact sets in $\mathcal{E}$. Here $\hat{\mathcal{E}}$ denotes the class of relatively compact sets in $\mathcal{E}$. We also let $\hat{\mathcal{E}}_\mu$ denote the sets in $\hat{\mathcal{E}}$ whose boundary has $\mu$-measure 0.

Let $\mathcal{C}_K^+$ denote the space of nonnegative functions on $\mathbb{E}$ with compact support. The integral of a functions $f \in \mathcal{C}_K^+$ with respect to a point process $\xi$ is the summation

$$\xi f \equiv \int_\mathbb{E} f(x)\, d\xi(x) \equiv \sum_k f(X_k).$$

Similarly, $\mu f \equiv \int_\mathbb{E} f(x)\, d\mu(x)$ for a measure $\mu$.

A sequence of point processes $\xi_n$ converges in distribution to the point process $\xi$, denoted by $\xi_n \xrightarrow{d} \xi$, if the probability measure of $\xi_n$ converges weakly to the probability measure of $\xi$. Each of the following statements is equivalent to $\xi_n \xrightarrow{d} \xi$: (a) $\xi_n f \xrightarrow{d} \xi f$, $f \in \mathcal{C}_K^+$. (b) $E[e^{-\xi_n f}] \to E[e^{-\xi f}]$, $f \in \mathcal{C}_K^+$. (c) $(\xi_n(B_1), \ldots, \xi_n(B_k)) \xrightarrow{d} (\xi(B_1), \ldots, \xi(B_k))$, $B_1, \ldots, B_k \in \hat{\mathcal{E}}_\mu$. When the limit $\xi$ is Poisson with mean measure $\mu$, statement (b) for the Laplace functional convergence is equivalent to

$$E[e^{-\xi_n f}] \to \exp\left\{-\int_\mathbb{E}(1 - e^{-f(x)})\mu(dx)\right\}.$$

We are now ready for the main result. Let $\xi_{n1}, \xi_{n2}, \ldots$, for $n \geq 1$, be a finite or infinite doubly indexed array of point processes on a Polish space $\mathbb{E}$ with Borel $\sigma$-field $\mathcal{E}$, so that $\xi_{nk}(B)$ is the number of points in $B \in \mathcal{E}$. Assume the processes are defined on a single probability space $(\Omega, \mathcal{F}, P)$, and there are $\sigma$-fields $\mathcal{F}_n$, $n \geq 1$, in $\mathcal{F}$ such that $\xi_{n1}, \xi_{n2}, \ldots$ are conditionally independent given $\mathcal{F}_n$. Assume that either

(1) $$\sum_k P\{\xi_{nk}(B) \geq 1 | \mathcal{F}_n\}^2 \xrightarrow{d} 0 \qquad \text{as } n \to \infty, \ B \in \hat{\mathcal{E}},$$

or that

(2) $$\sup_k P\{\xi_{nk}(B) \geq 1 | \mathcal{F}_n\} \xrightarrow{d} 0 \qquad \text{as } n \to \infty, \ B \in \hat{\mathcal{E}}.$$

Consider the sequence

$$N_n = \sum_k \xi_{nk}, \qquad n \geq 1, \tag{3}$$

and, for the infinite sequence case, assume that $N_n$ is finite on compact sets, so that it is a point process on $\mathbb{E}$. Let $N$ denote a Poisson process on $\mathbb{E}$ with mean measure $\mu$. The following result gives conditions under which $N_n \xrightarrow{d} N$.

THEOREM 1. *For the processes defined above,*

$$P\{N_n \in \cdot | \mathcal{F}_n\} \xrightarrow{d} P\{N \in \cdot\} \qquad as \ n \to \infty, \tag{4}$$

*if and only if the following conditions hold as $n \to \infty$:*

$$\sum_k P\{\xi_{nk}(B) \geq 2 | \mathcal{F}_n\} \xrightarrow{d} 0, \qquad B \in \hat{\mathcal{E}}, \tag{5}$$

$$\sum_k P\{\xi_{nk}(B) \geq 1 | \mathcal{F}_n\} \xrightarrow{d} \mu(B), \qquad B \in \hat{\mathcal{E}}_\mu. \tag{6}$$

*Moreover, these conditions imply $N_n \xrightarrow{d} N$, as $n \to \infty$.*

This result is a conditional-distribution variation of the classical theorem of Grigelionis [2] (which is Theorem 16.18 in [4]). In particular, he proved Theorem 1 for independent $\xi_{n1}, \xi_{n2}, \ldots$ that satisfy the uniformly null property (2), and all the conditional distributions are ordinary ones. In this case, his result states that (5) and (6) are necessary and sufficient for $N_n \xrightarrow{d} N$.

In addition to the standard uniformly null property (2), we included the alternative property (1), which we use later.

Our proof below of Theorem 1 uses ideas and techniques associated with the convergence of point processes (e.g., see Lemma 5.8 and the proofs of Theorems 5.7 and 16.18 in [4]). The following preliminary result plays a key role.

LEMMA 2. *Let $Y_{n1}, Y_{n2}, \ldots$, for $n \geq 1$, be a (possibly finite) doubly indexed array of random variables in $(0, c]$, where $c < 1$, and let $Y$ be a nonnegative random variable. Assume*

$$\sum_k Y_{nk}^2 \xrightarrow{d} 0 \quad or \quad \sup_k Y_{nk} \xrightarrow{d} 0 \qquad as \ n \to \infty. \tag{7}$$

*Then as $n \to \infty$,*

$$\prod_k (1 - Y_{nk}) \xrightarrow{d} e^{-Y} \quad if \ and \ only \ if \quad \sum_k Y_{nk} \xrightarrow{d} Y.$$

PROOF. The assertion is equivalent to

$$-\sum_k \log(1-Y_{nk}) \xrightarrow{d} Y \quad \text{if and only if} \quad \sum_k Y_{nk} \xrightarrow{d} Y. \tag{8}$$

To prove this equivalence, it suffices to show $D_n \xrightarrow{d} 0$, where

$$D_n = -\sum_k \log(1-Y_{nk}) - \sum_k Y_{nk} = \sum_k Y_{nk}^2 \left( \sum_{m=2}^{\infty} (Y_{nk}^{m-2}/m) \right).$$

Using $Y_{nk} \leq c$ and $\alpha_n \equiv \sup_k Y_{nk}$, we have

$$D_n \leq (1-c)^{-1} \sum_k Y_{nk}^2, \tag{9}$$

$$D_n \leq \frac{\alpha_n}{1-c} \sum_k Y_{nk} \leq \frac{-\alpha_n}{1-c} \sum_k \log(1-Y_{nk}). \tag{10}$$

Then $D_n \xrightarrow{d} 0$ follows from (9) under the supposition $\sum_k Y_{nk}^2 \xrightarrow{d} 0$. Also, under the alternative supposition $\alpha_n \xrightarrow{d} 0$, if one of the limit statements in (8) is true, then $D_n \xrightarrow{d} 0$ by (10), which proves the other limit statement in (8). □

PROOF OF THEOREM 1. We begin by determining a convenient expression for the statement $P\{N_n \in \cdot|\mathcal{F}_n\} \xrightarrow{d} P\{N \in \cdot\}$. This convergence is equivalent to the random Laplace functional convergence (e.g., see Theorems 5.3 and 16.16 of [4])

$$E[e^{-N_n f}|\mathcal{F}_n] \xrightarrow{d} E[e^{-Nf}], \qquad f \in \mathcal{C}_K^+. \tag{11}$$

Since $\xi_{n1}, \xi_{n2}, \ldots$ are conditionally independent given $\mathcal{F}_n$,

$$E[e^{-N_n f}|\mathcal{F}_n] = \prod_k E[e^{-\xi_{nk} f}|\mathcal{F}_n] = \prod_k (1-Y_{nk}),$$

where

$$Y_{nk} = E[1 - e^{-\xi_{nk} f}|\mathcal{F}_n].$$

Also, the Laplace functional of the Poisson process $N$ has the well-known form $E[e^{-Nf}] = e^{-\mu h}$, where $h(x) := 1 - e^{-f(x)}$. Then (11) is equivalent to

$$\prod_k (1-Y_{nk}) \xrightarrow{d} e^{-\mu h}, \qquad f \in \mathcal{C}_K^+. \tag{12}$$

Keep in mind that $Y_{nk}$ is a function of $f$.

We will complete the proof by applying Lemma 2 to establish that (5) and (6) are necessary and sufficient for (12). Clearly, $Y_{n1}, Y_{n2}, \ldots$ are in $(0, 1 - e^{-b}]$, where $b = \max_{x \in \mathbb{E}} f(x)$. Next, note that

$$(13) \quad 1 - e^{-\xi_{nk}f} \leq \mathbf{1}(\xi_{nk}(S_f) \geq 1),$$

where $S_f$ is the support of $f$. Then by assumption (1),

$$\sum_k Y_{nk}^2 \leq \sum_k P\{\xi_{nk}(S_f) \geq 1|\mathcal{F}_n\}^2 \xrightarrow{d} 0.$$

Also, by (13) and assumption (2),

$$\sup_k Y_{nk} \leq \sup_k P\{\xi_{nk}(S_f) \geq 1|\mathcal{F}_n\} \xrightarrow{d} 0.$$

In light of the last two limit statements, Lemma 2 says that (12) is equivalent to

$$(14) \quad \sum_k (1 - E[e^{-\xi_{nk}f}|\mathcal{F}_n]) = \sum_k Y_{nk} \xrightarrow{d} \mu h, \qquad f \in \mathcal{C}_K^+.$$

Therefore, it remains to show that (5) and (6) are necessary and sufficient for (14). We show sufficiency first.

We can write

$$(15) \quad \sum_k Y_{nk} = \sum_k E[(1 - e^{-\xi_{nk}f})\mathbf{1}(\xi_{nk}(S_f) = 1)|\mathcal{F}_n]$$
$$+ \sum_k E[(1 - e^{-\xi_{nk}f})\mathbf{1}(\xi_{nk}(S_f) \geq 2)|\mathcal{F}_n].$$

The last sum is bounded by $\sum_k P\{\xi_{nk}(S_f) \geq 2)|\mathcal{F}_n\}$ which converges in distribution to 0 by assumption (5). The first sum on the right-hand side in (15) equals $\sum_k E[\tilde{\xi}_{nk}h|\mathcal{F}_n] = \eta_n h$, where

$$\tilde{\xi}_{nk}(B) = \xi_{nk}(B \cap S_f)\mathbf{1}(\xi_{nk}(S_f) = 1),$$
$$\eta_n(B) = \sum_k E[\tilde{\xi}_{nk}(B)|\mathcal{F}_n] = \sum_k P\{\xi_{nk}(B \cap S_f) = 1|\mathcal{F}_n\}.$$

Now, assumptions (5) and (6) imply $\eta_n \xrightarrow{d} \mu$, and so $\eta_n h \xrightarrow{d} \mu h$. Applying the preceding observations to (15) proves that (5) and (6) are sufficient for (14).

Now suppose (14) is true and recall that (14), (12) and (11) are equivalent. Then from (11) and switching from Laplace transforms to generating functions, it follows that, for any fixed $B \in \mathcal{E}_\mu$,

$$E[s^{N_n(B)}|\mathcal{F}_n] \xrightarrow{d} E[s^{N(B)}], \qquad s \in [0,1].$$

This in turn implies that (12) and, hence, (14) also hold with $e^{-\xi_{nk}f}$ replaced by $s^{\xi_{nk}(B)}$. In particular, the version of (14) says for $s \in [0,1]$ that

$$(16) \qquad H_n(s) \equiv \sum_k [1 - E(s^{\xi_{nk}(B)}|\mathcal{F}_n)] \xrightarrow{d} (1-s)\mu(B).$$

Here we are using the fact that $E[e^{-\alpha X_n}] \to E[e^{-\alpha X}]$, $\alpha \in [0,\infty)$, implies $X_n \xrightarrow{d} X$, which implies $E[s^{X_n}] \to E[s^X]$, $s \in [0,1]$. Then (6) follows since

$$(17) \qquad \sum_k P\{\xi_{nk}(B) \geq 1|\mathcal{F}_n\} = H_n(0) \xrightarrow{d} \mu(B).$$

Next, note that

$$H_n(s) = \sum_k \left[1 - \sum_{m=0}^\infty s^m P\{\xi_{nk}(B) = m|\mathcal{F}_n\}\right]$$

$$= (1-s)H_n(0) + \sum_k \sum_{m=2}^\infty (s - s^m)P\{\xi_{nk}(B) = m|\mathcal{F}_n\}.$$

Then using this expression along with (16) and (17),

$$(s - s^2)\sum_k P\{\xi_{nk}(B) \geq 2|\mathcal{F}_n\} \leq H_n(s) - (1-s)H_n(0) \xrightarrow{d} 0.$$

These observations prove that (14) implies (5) and (6).

This completes the proof of the first assertion that (5) and (6) are necessary and sufficient for (4). Finally, (4) and the bounded convergence theorem yield $N_n \xrightarrow{d} N$, which proves the second assertion. $\square$

We will now consider the point processes defined above, where $\xi_{nk}$ are single-atom measures. Suppose

$$(18) \qquad N_n(B) = \sum_k U_{nk}\delta_{X_{nk}}(B), \qquad B \in \mathcal{E},$$

where $X_{n1}, X_{n2},\ldots$ are $\mathbb{E}$-valued random variables, and $U_{n1}, U_{n2},\ldots$ are random variables that take on values 0 or 1. Notice that (18) is the same as (3) with $\xi_{nk}(B) = U_{nk}\delta_{X_{nk}}(B)$. Assume that, for each $n \geq 1$, there is a $\sigma$-field $\mathcal{F}_n$ in $\mathcal{F}$ that contains $\sigma(X_{n1}, X_{n2},\ldots)$, and $U_{n1}, U_{n2},\ldots$ are conditionally independent given $\mathcal{F}_n$. Define

$$r_{nk} \equiv P\{U_{nk} = 1|\mathcal{F}_n\}, \qquad k \geq 1.$$

COROLLARY 3. *Suppose the $N_n$ defined in (18) are such that*

$$(19) \qquad \sum_k r_{nk}^2 \delta_{X_{nk}}(\cdot) \xrightarrow{d} 0 \quad \text{and} \quad \sum_k r_{nk}\delta_{X_{nk}}(\cdot) \xrightarrow{d} \mu, \qquad \text{as } n \to \infty,$$

for some measure $\mu$ that is finite on sets in $\hat{\mathcal{E}}$. Then $N_n \xrightarrow{d} N$, as $n \to \infty$, where $N$ is a Poisson process with mean measure $\mu$.

PROOF. Note that $P\{U_{nk}\delta_{X_{nk}}(B) = 1|\mathcal{F}_n\} = r_{nk}\delta_{X_{nk}}(B)$. Then the conditions in (19) imply (1) and (6), and (5) is trivially satisfied. Then $N_n \xrightarrow{d} N$ follows by Theorem 1. □

EXAMPLE 4. *Random variable convergence.* As a special case of (18), suppose $N_n = \sum_k U_{nk}$, where $U_{n1}, U_{n2}, \ldots$ are random variables that take on values 0 or 1. Assume that, for each $n \geq 1$, there is a $\sigma$-field $\mathcal{F}_n$ in $\mathcal{F}$ such that $U_{n1}, U_{n2}, \ldots$ are conditionally independent given $\mathcal{F}_n$. In addition, assume that as $n \to \infty$,

$$\sum_k P\{U_{nk} = 1|\mathcal{F}_n\}^2 \xrightarrow{d} 0 \quad \text{and} \quad \sum_k P\{U_{nk} = 1|\mathcal{F}_n\} \xrightarrow{d} \mu,$$

for some $\mu > 0$. Then by Corollary 3, the random variables $N_n$ converge in distribution to a Poisson random variable with mean $\mu$.

The next result is a special case of Corollary 3 that we use for the particle system model in the next section. Suppose

$$(20) \qquad N_n(B) = \sum_{k=1}^{n} U_{nk}\delta_{X_k}(B), \qquad B \in \mathcal{E},$$

where $X_1, X_2, \ldots$ are independent identically distributed random elements in $\mathbb{E}$ with distribution $F$, and $U_{n1}, \ldots, U_{nn}$ are random variables that take on values 0 or 1 and are conditionally independent given $\mathcal{F}_n \equiv \sigma(X_1, \ldots, X_n)$.

COROLLARY 5. *Suppose the point processes $N_n$ given by (20) satisfy*

$$(21) \qquad\qquad\qquad (r_{nk}, X_k) \stackrel{d}{=} (r_{n1}, X_1), \qquad 1 \leq k \leq n,$$

$$(22) \qquad \lim_{n \to \infty} n \int_B E[r_{n1}^2|X_1 = x]\,dF(x) = 0, \qquad\qquad B \in \hat{\mathcal{E}}.$$

*In addition, assume there is a function $r : \mathbb{E} \to [0, \infty)$ such that*

$$(23) \qquad \lim_{n \to \infty} \int_B E[|nr_{n1} - r(x)||X_1 = x]\,dF(x) = 0, \qquad B \in \hat{\mathcal{E}},$$

*and the measure $\mu(B) \equiv \int_B r(x)\,dF(x)$, $B \in \mathcal{E}$, is finite on compact sets. Then $N_n \xrightarrow{d} N$ as $n \to \infty$, where $N$ is a Poisson process with mean measure $\mu$.*

PROOF. The assertion will follow by Corollary 3 upon verifying (19). For any $f \in \mathcal{C}_K^+$, with support $S_f$, it follows by (21) and (22) that

$$E\left[\sum_{k=1}^n r_{nk}^2 f(X_k)\right] = E[nr_{n1}^2 f(X_1)]$$

$$\leq \max_x f(x) \int_{S_f} nE[r_{n1}^2 | X_1 = x] \, dF(x) \to 0.$$

This implies $\sum_{k=1}^n r_{nk}^2 f(X_k) \xrightarrow{d} 0$, which proves the first condition in (19).

To prove the second condition in (19), it suffices to show the $L_1$-convergence

(24) $$E\left|\sum_{k=1}^n r_{nk} f(X_k) - \mu f\right| \to 0, \qquad f \in \mathcal{C}_K^+.$$

By the classical $L_1$ law of large numbers, the independence of $X_1, X_2, \ldots$, and $E[r(X_1)f(X_1)] = \mu f$, we have

$$E\left|n^{-1}\sum_{k=1}^n r(X_k)f(X_k) - \mu f\right| \to 0.$$

Also, by (21) and (23),

$$E\left[\sum_{k=1}^n r_{nk} f(X_k) - n^{-1}\sum_{k=1}^n r(X_k)f(X_k)\right]$$

$$\leq \max_x f(x) E[|nr_{n1} - r(X_1)|] \to 0.$$

The last two limit statements imply (24). $\square$

**3. Particle system with attribute-dependent survival.** Suppose particles (representing microbes, molecules, customers, etc.) arrive to a system at finite times $0 \leq T_1 \leq T_2 \leq \cdots$. The number of particles that arrive in the time interval $[0, t]$ is $A(t) \equiv \sum_{n=1}^\infty \mathbf{1}(T_n \leq t)$; this is finite when $T_n \to \infty$ w.p.1, which we assume for convenience. No other conditions are imposed on the dependency among the $T_n$. Without loss of generality (as explained in Proposition 9), we assume the system is empty at time 0.

Each particle that arrives is labeled by a real-valued attribute (or rank) $x$. We interpret $x$ as the "location" of the particle in the attribute space, and refer to a particle with attribute $x$ as an "$x$-particle." We let $X_n$ denote the attribute of the particle that arrives at time $T_n$, and assume $X_n$ are independent, continuous random variables, independent of the arrival times, with a common distribution function $F$. Then the space of the particle attributes is the open interval

$$\mathbb{E} \equiv \{x \in \Re : 0 < F(x) < 1\}.$$

The particles in the system (i.e., in the space $\mathbb{E}$) are subject to deletion by future arrivals as follows. Whenever a new particle arrives, it considers for deletion all particles currently in the system that have strictly lower attributes. Specifically, an arriving $y$-particle considers all $x$-particles with $x < y$ for deletion, and then it deletes each such $x$-particle with probability $a(x)$, independently of everything else. The probability that an $x$-particle survives $\ell$ deletion attempts is therefore $(1 - a(x))^\ell$. From the view of an $x$-particle, the probability it is deleted by a $y$-particle is $\mathbf{1}(x < y)a(x)$. Also, the probability that an $x$-particle is considered for deletion by the next arrival is $\overline{F}(x) = 1 - F(x)$, and the probability it is actually deleted by that arrival is $a(x)\overline{F}(x)$.

We denote the *lifetime* of the $X_n$-particle by $L_n$. That is, $L_n = \ell$ means the $X_n$-particle dies (exits the system) at the $\ell$th deletion attempt. Under the preceding assumptions, the pairs $(X_n, L_n)$, $n \geq 1$ are independent, identically distributed, independent of the arrival times and

$$(25) \qquad P\{L_n > \ell | X_n\} = (1 - a(X_n))^\ell, \qquad \ell \geq 1.$$

The only other assumption we make is that the function $1/a(x)$ is bounded on compact sets in $\mathbb{E}$. This implies the measure $\mu$ defined below in (27) is finite on compact sets.

We will consider the point process $N_t$ on $\mathbb{E}$, where $N_t(B)$ denotes the number of particles with attributes in $B$ that are still in the system at time $t$. That is,

$$(26) \qquad N_t(B) = \sum_{k=1}^{A(t)} \mathbf{1}(L_k > Q_{tk})\delta_{X_k}(B), \qquad B \in \mathcal{E},$$

where $Q_{tk} = \sum_{j=k+1}^{A(t)} \mathbf{1}(X_j > X_k)$ is the quantity of attempted deletions of the $X_k$-particle by those $X_j$-particles that are above $X_k$ and arrive in the interval $[0, t]$ after $X_k$. Recall that the dynamics are such that the $X_k$-particle is alive at time $t$ if and only if its lifetime $L_k$ exceeds $Q_{tk}$.

The time-dependent evolution of the particle process $\{N_t : t \geq 0\}$ is highly dependent on the nature of the arrival process $A(t)$, which may have any probabilistic structure. However, when the arrival process is Poisson, then the particle process is a pure-jump Markov process. Specifically, the particle process changes state only at the Poisson arrival times $T_n$, and the embedded sequence of states $N_{T_n}$ at the arrival times is a Markov chain. That is, $N_t$ is a Markov chain subordinated to the Poisson process $A(t)$.

In this Markovian setting with the additional assumptions that the deletion probability $a(x) \equiv a$ is independent of the location $x$ and $F$ is a uniform distribution on $\mathbb{E} = (0, 1)$, Durrett and Limic [1] showed that $N_t$ converges in distribution to a Poisson process with mean measure $\mu(B) = \int_B \frac{1}{a(1-x)}\, dx$.

The approach they used to prove this is not applicable in the non-Markovian setting with location-dependent deletions.

We are now ready to apply the results above for sums of point processes to show the particle process $N_t$ converges in distribution to a Poisson process.

THEOREM 6. *Under the preceding conditions, $N_t \xrightarrow{d} N$ as $t \to \infty$, where $N$ is a Poisson process on $\mathbb{E}$ with mean measure*

$$(27) \qquad \mu(B) = \int_B \frac{1}{a(x)\overline{F}(x)} \, dF(x), \qquad B \in \mathcal{E}.$$

Intricate dependencies in the particle process cloud the fact that the Poisson limit of $N_t$ follows by using the framework in Corollary 5. Furthermore, the main property producing the Poisson limit is the *claim* (29) below about ranks of order statistics (but in terms of conditional distributions). Further insights underlying Poisson limits are in Section 5, where we discuss stationary distributions.

PROOF OF THEOREM 6. We begin with the key observation that, conditioned on $A(t) = n$, the process $N_t$ depends on time only through the "order" in which the $n$ particles arrive. Because of this observation, we can formulate the distribution of $N_t$ as follows.

Fix an $n \geq 1$, and let $(X_1, L_1), \ldots, (X_n, L_n)$ be i.i.d. versions of the random vectors above, but here the subscripts are not the order of particle arrivals. On the same probability space and independent of the $(X_k, L_k)$, let $\pi_{n1}, \ldots, \pi_{nn}$ be the random permutation of $1, \ldots, n$, where each permutation is equally likely with probability $1/n!$. We define the point process $\xi_n$ on $\mathbb{E}$ by

$$(28) \qquad \xi_n = \sum_{k=1}^n U_{nk} \delta_{X_k},$$

where $U_{nk} = \mathbf{1}(L_k > Q_{nk})$ and

$$Q_{nk} = \sum_{j=1}^n \mathbf{1}(X_j > X_k, \pi_{nj} > \pi_{nk}), \qquad 1 \leq k \leq n.$$

We interpret $X_1, \ldots, X_n$ as the attributes of the first $n$ particles to arrive, with $\pi_{n1}, \ldots, \pi_{nn}$ representing the order in which they arrive to the system; for example, $\pi_{nk} = 3$ means that the $X_k$-particle is the third one to arrive. Also, $Q_{nk}$ is the quantity of attempted deletions of the $X_k$-particle.

From the definition of the process $N_t$, it is clear that $\xi_n$ represents $N_t$ conditioned on $A(t) = n$. In other words,

$$P\{N_t \in \cdot\} = E[P_{A(t)}(\cdot)],$$

where $P_n(\cdot) = P\{\xi_n \in \cdot\}$, $n \geq 0$, is the distribution of $\xi_n$, and $\xi_0$ is the zero measure. Since $P_n$ does not depend on $t$ and $A(t) \to \infty$ w.p.1, the assertion $N_t \xrightarrow{d} N$ will follow upon showing that $P\{\xi_n \in \cdot\}$ converges weakly to $P\{N \in \cdot\}$, or equivalently, that $\xi_n \xrightarrow{d} N$.

We will prove this convergence by applying Corollary 5 to the process $\xi_n$ in (28). We begin by proving the following assertion.

CLAIM. *The $Q_{n1}, \ldots, Q_{nn}$ are conditionally independent given $\mathcal{F}_n = \sigma(X_1, \ldots, X_n)$, and*

$$(29) \quad P\{Q_{nk} = m | \mathcal{F}_n\} = 1/\nu_{nk}, \qquad 0 \leq m \leq \nu_{nk} - 1, 1 \leq k \leq n,$$

*where $\nu_{nk} = \sum_{j=1}^{n} \mathbf{1}(X_j \geq X_k)$.*

PROOF. It suffices to prove the claim for $\{X_1 > \cdots > X_n\}$. In this case, $\nu_{nk} = k$ and

$$Q_{nk} = \sum_{j=1}^{k-1} \mathbf{1}(\pi_{nj} > \pi_{nk}),$$

which is the number of the particles $X_1, \ldots, X_{k-1}$ that arrive after $X_k$. Since $\pi_{n1}, \ldots, \pi_{nn}$ is the equally likely permutation of $1, \ldots, n$, we have

$$P\{Q_{nk} = m_k | \mathcal{F}_n\} = 1/\nu_{nk}, \qquad 0 \leq m_k \leq \nu_{nk} - 1,$$

(30)
$$P\left\{\bigcap_{k=1}^{n} \{Q_{nk} = m_k\} | \mathcal{F}_n\right\} = 1/n! = \prod_{k=1}^{n} 1/\nu_{nk}.$$

The first equality of the second line above follows because, conditioned on $\mathcal{F}_n$, each of the $n!$ equally likely permutations of arrival orderings corresponds to exactly one distinct assignment of $Q_{nk}$ values. This proves the claim. □

Next, we derive a simple expression for $r_{nk} \equiv P\{U_{nk} = 1 | \mathcal{F}_n\}$. Using (25) and $\bar{a}(x) = 1 - a(x)$, we have

$$r_{nk} \equiv E[P\{U_{nk} = 1 | Q_{nk}, \mathcal{F}_n\} | \mathcal{F}_n] = E[\bar{a}(X_{nk})^{Q_{nk}} | \mathcal{F}_n].$$

Then from (30),

$$(31) \quad r_{nk} = \frac{1}{\nu_{nk}} \sum_{m=0}^{\nu_{nk}-1} \bar{a}(X_{nk})^m = \frac{[1 - \bar{a}(X_{nk})^{\nu_{nk}}]}{\nu_{nk} a(X_{nk})}.$$

We are now ready to verify the assumptions in Corollary 5. By assumption, $X_{nk} = X_k$, $1 \leq k \leq n$, are independent and identically distributed with distribution $F$, and $U_{n1}, \ldots, U_{nn}$ are conditionally independent given $\mathcal{F}_n$ by the claim we proved above. Also, $(r_{nk}, X_k) \stackrel{d}{=} (r_{n1}, X_1)$, $1 \leq k \leq n$, follows from (31) and the fact that $\nu_{nk} = \sum_{j=1}^{n} \mathbf{1}(X_j \geq X_k)$ is a function of $X_k$ and the set $\{X_1, \ldots, X_n\} \setminus \{X_k\}$.

Our next step is to prove

$$(32) \qquad \lim_{n \to \infty} n \int_B E[r_{n1}^2 | X_{n1} = x] \, dF(x) = 0, \qquad B \in \hat{\mathcal{E}}.$$

First note that

$$(33) \qquad P\{\nu_{n1} = m | X_1 = x\} = P\{S_{n-1} + 1 = m\},$$

where $S_{n-1}$ is a binomial random variable with parameters $n-1$ and $\overline{F}(x)$. Then

$$(34) \qquad nE[r_{n1}^2 | X_{n1} = x] = E\left[\frac{n[1 - \overline{a}(x)^{S_{n-1}+1}]^2}{(S_{n-1}+1)^2}\right].$$

To determine the limiting behavior of this expression, we will use the identity

$$(35) \qquad E[f(S_n)] = \frac{1}{p(n+1)} E[S_{n+1} f(S_{n+1} - 1) \mathbf{1}(S_{n+1} \geq 1)],$$

where $f : \{0, \ldots, n\} \to \Re$, and $S_n$ is a binomial random variable with parameters $n$ and $p$. This follows by writing out the expectations.

Applying (35) to (34) yields

$$nE[r_{n1}^2 | X_{n1} = x] = E\left[\frac{[1 - \overline{a}(x)^{S_n}]^2 \mathbf{1}(S_n \geq 1)}{S_n \overline{F}(x)}\right].$$

The last expression in brackets is bounded by $\mathbf{1}(S_n \geq 1)/(S_n \overline{F}(b))$, where $b = \sup\{x : x \in B\}$ and it converges to 0 w.p.1 since $\overline{F}(b) > 0$ and $S_n \to \infty$ w.p.1. Thus, (22) is true by the bounded convergence theorem.

Next, consider the function

$$r(x) \equiv \frac{1}{a(x)\overline{F}(x)}.$$

This is bounded on compact sets in $\mathbb{E}$ since $1/a(x)$ is, and so the measure $\mu$ in (27) is bounded on compact sets. It remains to prove (23) in Corollary 5, which, in light of (31) and (33), is

$$(36) \qquad \lim_{n \to \infty} \int_B E\left|\frac{n[1 - \overline{a}(x)^{S_{n-1}+1}]}{a(x)(S_{n-1}+1)} - r(x)\right| dF(x) = 0, \qquad B \in \hat{\mathcal{E}}.$$

Using (35), we see that the expectation in this integral equals

$$(37) \qquad r(x) E\left[\left|1 - \overline{a}(x)^{S_n} - \frac{S_n}{n\overline{F}(x)}\right| \mathbf{1}(S_n \geq 1)\right].$$

The expression inside the expectation is bounded by $1 + 1/\overline{F}(b)$, where $b = \sup\{x : x \in B\}$ and it converges to 0 because of the classical strong law of large numbers $n^{-1} S_n \to \overline{F}(x)$ w.p.1. Then the expectation in (37) converges to 0 by the bounded convergence theorem. Furthermore, since $r(x)$ is bounded on $B$, the dominated convergence theorem justifies (36).

This completes the verification of conditions in Corollary 5, which yields $\xi_n \xrightarrow{d} N$ as $n \to \infty$, and hence, $N_t \xrightarrow{d} N$ as $t \to \infty$. □

**4. Sojourn and departure times of particles.** This section describes sojourn times of particles in the system defined above. The sojourn times have tractable distributions when the arrival times form a renewal or Poisson process. We will also discuss some features of the departure times of particles.

Consider the particle process $\{N_t : t \geq 0\}$ described in the preceding section, where the initial state $N_0$ need not be 0, and $N_t$ may be stationary or nonstationary (in the time parameter). We will consider the sojourn time $W(x)$ of a "typical" $x$-particle that arrives at a time $t$. This sojourn time is "conditioned" on having an arrival at a fixed time $t$ with attribute $x$. Conditioning on such events (which may have probability 0) requires the use of a family of Palm probabilities $P_t$ ($t \geq 0$) for nonstationary systems as described in [3]. When a system is stationary, these Palm probabilities are independent of $t$ and equal to a single probability measure. Since a formal definition of these probabilities is lengthy, we will not include it; the following discussion can be understood without a working knowledge of Palm probabilities.

For the next result, let $P_t$ denote the Palm probability of the entire system conditioned that a particle arrives at time $t$. In addition to $W(x)$, we will consider the *spatial* sojourn time $W(B)$ of a typical particle whose attribute is in the set $B \in \mathcal{E}$.

PROPOSITION 7 (Sojourn times). *Suppose the arrival process $A(t)$ for the particle system is a renewal process. Then the Palm distribution of the sojourn time $W(B)$ for a particle entering $B$ at time $t$ is*

$$(38) \quad P_t\{W(B) > w\} = \int_B E[(1 - a(x)\overline{F}(x))^{A(w)}] \, dF(x), \qquad w \geq 0.$$

*In case $A(t)$ is a Poisson process with rate $\lambda$,*

$$(39) \quad P_t\{W(B) > w\} = \int_B e^{-\lambda a(x)\overline{F}(x)} \, dF(x),$$

*and hence, the sojourn time $W(x)$ for an $x$-particle entering at time $t$ has an exponential distribution with rate $\lambda a(x) \overline{F}(x)$.*

PROOF. By the definitions of the random variables involved,

$$P_t\{W(x) > w\} = P\{A(w) < \nu(x)\}, \tag{40}$$

where $\nu(x)$ is the number of new arrivals until the $x$-particle is deleted. Under the rule for particle deletions, an $x$-particle is deleted by an arbitrary arrival with probability $a(x)\overline{F}(x)$, independently of everything else. Then $\nu(x)$ has a geometric distribution with $P\{\nu(x) > n\} = (1 - a(x)\overline{F}(x))^n$, and $\nu(x)$ is independent of $A(w)$. Applying these properties to (40), we have

$$P_t\{W(x) > w\} = E[(1 - a(x)\overline{F}(x))^{A(w)}]. \tag{41}$$

Then conditioning on $x$ being in the set $B$ yields (38).

In case $A(t)$ is a Poisson process with rate $\lambda$, we know that $E[s^{A(t)}] = e^{-\lambda(1-s)}$, and so (41) reduces to (39). □

In the preceding result, the sojourn-time distribution of a particle is independent of the time at which the particle enters the system. This independence of time is because the sojourn of a particle does not depend on any of the system data from the past. It is a rare instance in which a certain probability for a nonstationary system is the same as the analogous probability for the stationary version of the system. Typically, Palm probabilities for a nonstationary system are not equal to those for its stationary version. For example, in most input–output systems or queueing networks, customer sojourn times depend highly upon quantities such as how many particles (or customers) are in the system upon arrival, as well as the total waiting times of the particles (customers) present.

REMARK 8. Even though sojourn times of particles in the setting of Proposition 7 are identically distributed, they are not independent. This follows because the simultaneous deletion of particles triggered by arrivals introduce dependencies among the sojourns.

PROPOSITION 9. *The Poisson convergence $N_t \xrightarrow{d} N$ in Theorem 6 is also true when the initial population process $N_0$ is any nonzero point process on $\mathbb{E}$.*

PROOF. Because the particles do not move about in $\mathbb{E}$, and the distribution of a point process is determined by its values on compact sets, the assertion will follow upon showing that, for any compact set $C \in \mathcal{E}$, the maximum sojourn time $\tau(C)$ of the $N_0(C)$ particles in $C$ at time 0 is a finite random variable.

Arguing as in the proof of Proposition 7 and conditioning on $N_0(\cdot)$ and $A(w)$, we have

$$P\{\tau(C) < w\} \geq E[(1 - (1-b)^{A(w)})^{N_0(C)}],$$

where $X_k$ are the locations of the $N_0(C)$ particles and $b \equiv \inf_{x \in C} a(x)\overline{F}(x)$. Because $N_0$ is a point process, $N_0(C)$ is finite. Also, $A(w) \to \infty$ w.p.1 as $w \to \infty$ by assumption. Thus,

$$P\{\tau(C) < \infty\} \geq \lim_{w \to \infty} E[(1-(1-b)^{A(w)})^{N_0(C)}] = 1$$

by dominated convergence. This completes the proof. $\square$

Another aspect of interest for the particle system is its departure process. The departure processes for Jackson and related queueing networks in equilibrium are space–time Poisson processes (e.g., see [6]). The departure times in the particle system are a little more complicated. The particles depart in batches at the arrival times of particles, and each batch is represented by a point process on $\mathbb{E}$ describing the locations where the particles depart.

Specifically, the departures are described by the space–time point process $M$ on $\Re_+ \times \mathbb{E}$ of the form

$$M(I \times B) = \sum_n \mathbf{1}(T_n \in I) M_n(B),$$

where $M_n(\cdot)$ is the point process of batch departures at time $T_n$. Under the assumptions for particle deletions, $M$ is a space-time $\mathcal{F}_t$-Markov process with respect to $\mathcal{F}_t \equiv \sigma(N_s(\cdot) : s \leq t)$ [which contains $\sigma(A(s) : s \leq t)$]. In particular, the increments of $M$ occur at the Poisson arrival times with rate $\lambda$, and the increments $M_1, M_2, \ldots$ at these times form a $\mathcal{F}_{T_n}$-Markov chain, and $M_n$ conditioned on $\mathcal{F}_{T_{n-1}}$ is a point process with intensity $a(x)\mathbf{1}(X_n > x)N_{T_{n-1}}(dx)$.

Even when the system is in equilibrium, the increments $M_n$ of the departure process $M$ are dependent. However, a "typical" batch departure process $M_n$ is a Poisson process as described below.

PROPOSITION 10. *Suppose the arrival process for the particle system is a Poisson process and the system is in equilibrium, and let $M_0(\cdot)$ denote the batch departure process conditioned on an arrival to the system at time 0 (as defined by the Palm probability). Then the spatial batch-departure process $M_0$ under the Palm probability is a Poisson process on $\mathbb{E}$ with intensity measure $a(x)\overline{F}(x)\,d\mu(x)$.*

PROOF. The Poisson limiting distribution for $N_t$ described by Theorem 6 is its stationary distribution, since this limit does not depend on the distribution of $N_0$ by Proposition 9. So the state of the process in equilibrium is given by a Poisson process $N$ with mean measure $\mu$ as in Theorem 6. Then under the Palm probability of the process conditioned on an arrival at time 0, the spatial batch-departure process $M_0$ is a conditional $p(x)$-thinning of $N$, where $p(x) = 1 - a(x)\overline{F}(x)$, which is the probability that an $x$-particle is retained in the system. Consequently, $M_0$ is a Poisson process with intensity measure $a(x)\overline{F}(x)\,d\mu(x)$. $\square$

**5. Stationary distribution of Markovian particle systems.** In the particle system we are studying, the probability that an $x$-particle is deleted by an arriving $y$-particle has the very special form $p(x,y) \equiv \mathbf{1}(x<y)a(x)$. A natural question is "Are there more general deletion probabilities $p(x,y)$ under which $N_t$ has a Poisson limit?" Also, the attribute space $\mathbb{E}$ for the particles is an interval in $\Re$ that is linearly ordered. How important is this linear order for Poisson limits?

We will give some insight into these issues with the following results. Consider a particle system somewhat like the one we have studied that evolves as follows. Particles arrive to the system according to a Poisson process $A(t)$ with rate $\lambda$, and each particle has an attribute in a Polish space $\mathbb{E}$. The attributes are i.i.d. with distribution $F$. Whenever a new particle arrives with an attribute $y$, each particle in the system is subject to being deleted; $p(x,y)$ is the probability that an $x$ particle is deleted by the arriving $y$-particle, independently of everything else. The only assumptions we make on the probabilities $p(x,y)$ is that $\int_{\mathbb{E}} p(x,y)\,dF(y) > 0$ $F$-a.e. and that the measure

$$(42) \qquad \mu(B) = \int_B \frac{1}{\int_{\mathbb{E}} p(x,y)\,dF(y)}\,dF(x), \qquad B \in \mathcal{E},$$

is finite on compact sets.

Let $N_t(B)$ denote the number of particles with attribute in $B$ at time $t$. As above, because of the nature of the deletions and the Poisson arrivals, $\{N_t : t \geq 0\}$ is a Markov jump process. Here is a characterization of the mean measure of a point process whose distribution is the stationary distribution of $N_t$.

PROPOSITION 11. *Suppose the stationary distribution of the Markov process $\{N_t : t \geq 0\}$ exists and it is the distribution of a point process $N$ on $\mathbb{E}$ whose mean measure $\mu$ is finite on compact sets. Then $\mu$ is given by (42).*

PROOF. Suppose the Markov process $N_t$ is stationary, which is equivalent to $N_0 \stackrel{d}{=} N$. Since $N_t$ is a Markov chain subordinated to a Poisson process, its stationary distribution is the same as that of its embedded Markov chain $\tilde{N}_n \equiv N_{T_n}$. Therefore, $\tilde{N}_1 \stackrel{d}{=} N$.

Now, $\tilde{N}_1 = \xi + \delta_{X_1}$, where $X_1$ is the attribute of the particle that arrives at time $T_1$, and $\xi$ is the point process consisting of those points of $N_0$ that survive the deletions at time $T_1$ triggered by the arriving particle. Consequently,

$$(43) \qquad E[N(B)] = E[\tilde{N}_1(B)] = E[\xi(B)] + P\{X_1 \in B\}, \qquad B \in \mathcal{E}.$$

Under the assumptions on particle deletions, we know that when $X_1 = x_1$, the process $\xi$ is a $p(x, x_1)$-thinning of $N_0 \stackrel{d}{=} N$. Then conditioning on $X_1$,

$$E[\xi(B)] = E\left[\int_B [1 - p(x, X_1)] \, dN(x)\right]$$
$$= \mu(B) - \int_B E[p(x, X_1)] \, d\mu(x).$$

Substituting this in (43) yields

$$\int_B E[p(x, X_1)] \, d\mu(x) = P\{X_1 \in B\}.$$

Thus, $\mu$ is given by (42). $\square$

We will now consider conditions under which the stationary distribution of the Markov process $N_t$ is that of a (spatial) Poisson process on $\mathbb{E}$. We already have one example.

REMARK 12. Suppose $N_t$ is the particle process as in Section 3 with a Poisson arrival process with rate $\lambda$. The stationary distribution of $N_t$ is that of a Poisson process $N$ with mean measure given by (27). This follows since $N_t \stackrel{d}{\to} N$ by Theorem 6, and this limit is independent of the distribution of $N_0$ by Proposition 9.

PROPOSITION 13. *The stationary distribution of the Markov process $\{N_t : t \geq 0\}$ is that of a Poisson process on $\mathbb{E}$ with mean measure (42) if and only if*

$$(44) \qquad \int_{\mathbb{E}} e^{\int_{\mathbb{E}} (1 - e^{-f(x)}) h(x, x_1) \, dF(x)} e^{-f(x_1)} \, dF(x_1) = 1, \qquad f \in \mathcal{C}_K^+,$$

*where $h(x, x_1) = p(x, x_1) / \int_{\mathbb{E}} p(x, y) \, dF(y)$.*

PROOF. Suppose $N_0$ is a Poisson process on $\mathbb{E}$ with mean measure $\mu$ given by (42). Then the stationary distribution of the Markov process $N_t$ is that of the Poisson process $N_0$ if and only if $\tilde{N}_1 \stackrel{d}{=} N_0$. The latter, in terms of Laplace functionals, is equivalent to

$$(45) \qquad E[e^{-\tilde{N}_1 f}] = E[e^{-N_0 f}], \qquad f \in \mathcal{C}_K^+.$$

Then to prove the proposition, it suffices to show (45) is equivalent to (44). Since $N_0$ is a Poisson process,

$$E[e^{-N_0 f}] = e^{-\int_{\Re} (1 - e^{-f(x)}) \, d\mu(x)}, \qquad f \in \mathcal{C}_K^+.$$

Using the representation $\tilde{N}_1 = \xi + \delta_{X_1}$, from the preceding proof, and conditioning on $X_1$, we have

$$E[e^{-\tilde{N}_1 f}] = \int_{\mathbb{E}} E[e^{-\xi f}|X_1 = x_1] e^{-f(x_1)} F(dx_1). \tag{46}$$

By the assumptions on particle deletions, we know that when $X_1 = x_1$, the process $\xi$ is a $p(x, x_1)$-thinning of the Poisson process $N_0$. Then the well-known theorem on thinning of Poisson processes implies that $\xi$, conditioned on $X_1 = x_1$, is a Poisson process with mean measure $[1 - p(x, x_1)] d\mu(x)$. Therefore,

$$\log E[e^{-\xi f}|X_1 = x_1] = -\int_{\mathbb{E}} (1 - e^{-f(x)})[1 - p(x, x_1)] d\mu(x)$$
$$= \log E[e^{-N_0 f}] + g(x_1),$$

where

$$g(x_1) = \int_{\mathbb{E}} (1 - e^{-f(x)}) p(x, x_1) d\mu(x).$$

Substituting this expression in (46), we obtain

$$E[e^{-\tilde{N}_1 f}] = E[e^{-N_0 f}] \int_{\mathbb{E}} e^{-[f(x_1) - g(x_1)]} F(dx_1). \tag{47}$$

Therefore, (45) is equivalent to

$$\int_{\mathbb{E}} e^{-[f(x_1) - g(x_1)]} F(dx_1) = 1, \qquad f \in \mathcal{C}_K^+.$$

But this expression, with $d\mu(x) = dF(x)/\int_{\mathbb{E}} p(x, y) dF(y)$ from (42), is equivalent to (44). $\square$

EXAMPLE 14. The particle process in Section 3 satisfies condition (44). Hence, the stationary distribution of $N_t$ is that of a Poisson process, which we had already confirmed in Remark 12. This observation follows, since in this case,

$$p(x, y) = \mathbf{1}(x < y) a(x), \qquad h(x, x_1) = \mathbf{1}(x < x_1)/\overline{F}(x),$$

and (44) reduces to

$$\int_{\Re} e^{-\int_{-\infty}^{x_1} (e^{-f(x)}/\overline{F}(x)) dF(x)} \frac{e^{-f(x_1)}}{\overline{F}(x_1)} dF(x_1) = 1. \tag{48}$$

Here one uses $\int_{-\infty}^{x_1} \frac{dF(x)}{\overline{F}(x)} = -\log \overline{F}(x_1)$. Then clearly (48), with a change of variable, equals $\int_0^{\infty} e^{-u} du = 1$, which is true.

EXAMPLE 15. *Counterexample.* Consider a particle system that evolves as the process described in Section 2, except that an arrival to $y$ will delete an $x$-particle with probability $p(x,y) = a(x)$, independent of $y$. In this case $h(x, x_1) \equiv 1$ and (44) reduces to $\int_{\mathbb{E}} e^{-f(x)} dF(x) = 1$, which is not true for all $f \in \mathcal{C}_K^+$. Hence, the stationary distribution of $N_t$ cannot be that of a Poisson process.

**Acknowledgments.** The authors thank Brian Fralix for his insightful proofreading that led to significant improvements in Section 2. We also thank the reviewers for their time, comments and suggestions.

DEPARTMENT OF MATHEMATICS
AUSTIN PEAY STATE UNIVERSITY
CLARKSVILLE, TENNESSEE 37044
USA
E-MAIL: jonesmatt@apsu.edu

SCHOOL OF INDUSTRIAL AND SYSTEMS ENGINEERING
GEORGIA INSTITUTE OF TECHNOLOGY
ATLANTA, GEORGIA 30332
USA
E-MAIL: rserfozo@isye.gatech.edu